\documentclass[12pt]{article}
\usepackage[english]{babel}
\usepackage{amsmath}
\usepackage{amssymb}
\usepackage{amsfonts}
\usepackage{amsthm}
\usepackage{graphics,graphicx}

\theoremstyle{plain}

\newtheorem{theorem}{Theorem}
\newtheorem*{proposition*}{Proposition}

\newtheorem{lemma}{Lemma}

\newtheorem*{corollary*}{Corollary}
\newtheorem*{Corollary*}{Important corollary}
\newtheorem{remark}{Remark}

\theoremstyle{definition}

\theoremstyle{remark}

\newtheorem*{remark*}{Remark}

\newtheorem*{example*}{Example}

\newtheorem*{examples*}{Examples}

\numberwithin{equation}{section}
\numberwithin{lemma}{section}
\numberwithin{theorem}{section}
\numberwithin{hypothesis}{section}
\numberwithin{definition}{section}
\numberwithin{example}{section}
\numberwithin{corollary}{section}
\numberwithin{remark}{section}
\begin{document}

\begin{center}
{\large\bf An asymptotic method of factorization \\
of a class of matrix-functions}

\vspace{3mm}
{\bf Gennady MISHURIS$^1$, Sergei ROGOSIN$^{1, 2}$}

\vspace{3mm}
{\footnotesize\it $^1$Aberystwyth University, Penglais, SY23 2BZ Aberystwyth, UK;

\vspace{1mm}
e-mail: ggm@ber.ac.uk; ser14@aber.ac.uk}

\vspace{3mm}
{\footnotesize\it $^2$Belarusian State University, Nezavisimosti Ave., 4, 220030 Minsk, Belarus;

\vspace{1mm}
e-mail: rogosin@bsu.by}

\end{center}

{\footnotesize {\bf Abstract.}
A novel method of asymptotic factorization of $n \times n$ matrix functions is proposed.
Considered class of matrices is motivated by certain problems originated in the elasticity theory.
An example is constructed to illustrate effectiveness of the proposed procedure. Further applications of the method is
discussed.}



\section{Introduction}

We are considering here the problem of factorization of continuous
matrix functions of the real variable. It means the representation of a given invertible square matrix
$G \in\left( {\mathcal C}({\mathbb R})\right)^{n\times n}$ in the following form
\begin{equation}
\label{fact}
G(x) = G^{-}(x) \Lambda(x)  G^{+}(x),
\end{equation}
where continuous invertible matrices $G^{-}(x)$, $G^{+}(x)$ possess an analytic continuation in the lower
$\Pi^{-} = \{z = x + i y: {\mathrm{Im}}\, < 0\}$ and upper $\Pi^{+} = \{z = x + i y: {\mathrm{Im}}\, > 0\}$
half-planes, respectively, and
\begin{equation}
\label{fact1}
\Lambda(x) = {\mathrm{diag}}\, \left(\left(\frac{x - i}{x + i}\right)^{\kappa_1}, \ldots, \left(\frac{x - i}{x + i}\right)^{\kappa_n}\right),\; \kappa_1, \ldots, \kappa_n\in {\mathbb Z}.
\end{equation}
The representation (\ref{fact}) is called {\it right (continuous or standard) factorization} and can be considered for any oriented curve $\Gamma$ of certain classes which divides complex plane into two domain $D^{-}$ and $D^{+}$ with changing of diagonal entries in $\Lambda(x)$ for $\left(\frac{x - t^{+}}{x - t^{-}}\right)^{\kappa_j}$, $t^{\mp}\in D^{\mp}$, or for $x^{\kappa_j}$. Similar representation
$$
G(x) = G^{+}(x) \Lambda(x)  G^{-}(x)
$$
is called {\it left (continuous or standard) factorization}. If right- (left-) factorization exists then the integer numbers $\kappa_1, \ldots, \kappa_n$, called {\it partial indices}, are determined uniquely up to the order. The factors $G^{-}$, $G^{+}$ are not unique. Relations of pairs of factors are described, e.g., in \cite{LitSpi87}. In particular, there exist constant transformations of factors, such that $\kappa_1 \geq \ldots \geq \kappa_n$. The right- (left-) factorization is called {\it canonical factorization} if all partial indices are equal to 0, i.e. $\kappa_1 = \ldots = \kappa_n = 0$.

Factorization of matrix functions was first studied in relation to the vector-matrix Riemann (or Riemann-Hilbert) boundary value problem (see \cite{Gak77}). The later was formulated by Riemann in his work on construction of complex differential equations with algebraic coefficients having a prescribed monodromy group (see, e.g., \cite{EhrSpi01}). By using Cauchy type integral's method the  vector-matrix Riemann boundary value problem was reduced in \cite{Mus68}, \cite{Vek67}  to a system of the Fredholm integral equations.  A part of the theory of the factorization problem is based on the study of such systems (see also \cite{Gak52}) though this approach does not answer, in particular, on the questions when it is possible to get factorization, how to construct factors and how to determine partial indices.

Among other sources of interest to the factorization problems one can point out the vector-valued Wiener-Hopf equations on a half-line (see \cite{GohKre58}, \cite{Spe85}) and their discrete analogous, namely the block Toeplitz equations (see, e.g., \cite{BoeSil90},  \cite{GohFel74}).  The developed technique found several applications in diffraction theory, fracture mechanics, geophysics, financial mathematics etc. (see a brief description given, e.g., in \cite{VeiAbr07} and references therein).

Theoretical background for the study of the matrix factorization and its numerous generalizations is presented in \cite{BoeSpi13}, \cite{EhrSpi01}, \cite{LitSpi87}, \cite{Vek67} (see also \cite{LawAbr07}).

The theory of the factorization is more or less complete (see \cite{BoeSpi13}), but the above mentioned constructive questions about existence, factors and partial indices (which are very important for practical applications) have been answered only in a number of special cases.  Among them one can mention rational matrix functions (see, e.g., \cite{Gak52}), functional commutative matrix functions (those satisfying $G(t) G(s) = G(s) G(t), \forall t, s\in \Gamma$, see \cite{Che56a}), upper- (lower-) triangular matrices with factorizable diagonal elements (see \cite{Che56}, \cite{FelGohKru94}), certain classes of meromorphic matrix functions (see \cite{Adu92}, \cite{AmiKam07}, \cite{Kiy12}), special cases of $2\times 2$  Daniele-Khrapkov matrix functions (with a small degree of deviator polynomial) (see \cite{Dan78}, \cite{EhrSpe02}, \cite{Khr71a}, \cite{Khr71b}), special cases of $2\times 2$ matrix functions with three rationally independent entries (see \cite{Akt92}, \cite{EhrSpe02} and references therein), special cases of $n\times n$ generalization of  the Daniele-Khrapkov matrix functions (see, e.g., \cite{CaSaMa01}, \cite{Jon84a}, \cite{Raw85}, \cite{VeiAbr07}), special classes of matrices possessing certain symmetry property (see \cite{JanLagEph99}, \cite{Vor11}  and references therein).

In this paper we propose an asymptotic method of construction of factors for a special class of  $n\times n$ nonrational matrix functions. To the best of authors' knowledge, this class does not coincide with any of the above mentioned classes. This class contains matrix functions appeared at the study of certain problems in fracture mechanics related to perturbation of the crack propagation (\cite{MishMovMov10}, \cite{PiMiMo07} -- \cite{PiMiMo12}).
The main idea of our study is to reduce determination of factors at each step of approximation to the solution of so-called vector-matrix jump boundary value problem.
The paper is organized as follows. In Sec. 2 we introduce necessary notation and formulate the problem.
Constructive algorithm is presented in Sec. 3. We also find interesting to present here the realization of the algorithm in a special case of matrices of practical importance.  The method is illustrated by an example given in Sec. 4. We conclude our study by showing the quality of the factorization approximation by restricting ourselves only to the first asymptotic term and discuss a role of the chosen small parameter.

\section{A class of matrices. Problem formulation}

Let us introduce the following class of invertible continuous $n\times n, n \geq 2$, matrix-functions ${\mathcal G} K_n$
depending on a real parameter $\varphi\in {\mathbb R}$,
satisfying the following conditions:

\noindent (1) $G_{\varphi}\in \left({\mathcal C}({\mathbb R})\right)^{n\times n}$ belongs to ${\mathcal G} K_n$ if it can be represented in the form
\begin{equation}
\label{cond1}
G_{\varphi} = R_{\varphi} F R_{\varphi}^{-1},
\end{equation}
where bounded locally H\"older-continuous on ${\mathbb R}$ (in general non-rational) invertible matrix $R_{\varphi}(x)$ is such that

\noindent (2)
\begin{equation}
\label{cond2}
R_{0} = R_{\varphi}\bigl|_{\varphi=0} = I,
\end{equation}

\noindent (3) matrix function $F$ does not depend on parameter $\varphi$, has H\"older-continuous entries $f_{kl}$
on the extended real line $\overline{\mathbb R}$, i.e. $\forall k,l = 1, \ldots, n,$
\begin{equation}
\label{cond3}
\left|f_{kl}(x_1) - f_{kl}(x_2)\right| \leq C \left|\frac{1}{x_1 + i} - \frac{1}{x_2 + i}\right|^{\mu},\; \forall x_1, x_2\in \overline{\mathbb R}, \;
0 < \mu < 1,
\end{equation}
 satisfies the following
asymptotic estimate at infinity

\noindent (4)
\begin{equation}
\label{cond4}
F(x) \rightarrow I,\;\;\; |x| \rightarrow \infty,
\end{equation}

\noindent (5) $F$ admits a right canonical factorization, i.e.
\begin{equation}
\label{cond5}
F(x) = F^{-}(x) F^{+}(x),
\end{equation}
where H\"older-continuous on $\overline{\mathbb R}$ matrix-functions $F^{-}(x), F^{+}(x)$ possess an analytic continuation in the lower
${\Pi}^{-}$ and the upper ${\Pi}^{+}$ half-plane, respectively.

The matrices of the following form constitute a simple subclass of class ${\mathcal G} K_2$:
\begin{equation}
\label{example}
G_{\varphi}(x) = \left(
\begin{array}{lr}
p(x) & q(x) e^{i \varphi x} \\
 q(x) e^{-i \varphi x} & p(x)
 \end{array}
 \right).
\end{equation}
This function appears after Fourier transforms of the Wiener-Hopf equation describing a problem of fracture mechanics.

In this particular case
\[
R_\varphi(x)=
\left(
\begin{array}{lr}
e^{-i\varphi x/2} & 0 \\
0 & e^{i\varphi x/2}
 \end{array}
 \right).
\]

We note that the matrix-functions of this type do not belongs to any known class
of matrix-functions which admit explicit factorization.

It can be seen below that in our algorithm the factors retain one of the properties of the matrix-function
$G_{\varphi}$, namely,
\begin{equation}
\label{infinity_unit}
G_{\varphi}^{-}(z), G_{\varphi}^{+}(z) \rightarrow I,\; {\mathrm{as}}\; z \rightarrow \infty, \; \mp {\mathrm{Im}}\, z > 0.
\end{equation}

\section{An algorithm}

\subsection{General construction}
\label{general}

By assumptions any matrix $G_{\varphi}(x)\in {\mathcal G} K_n$ can be written in the form
\begin{equation}
\label{constr1}
G_{\varphi}(x) = F^{-}(x) \left(F^{-}(x)\right)^{-1} G_{\varphi} \left(F^{+}(x)\right)^{-1} F^{+}(x) = F^{-}(x) G_{1, \varphi} F^{+}(x),
\end{equation}
where $F^{-}(x), F^{+}(x)$ are component of canonical factorization of the corresponding matrix $F$, and (see (\ref{cond2}), (\ref{cond4}))  the matrix $G_{1, \varphi}(x)$ is represented in the form
$$
G_{1, \varphi}(x) = \left(F^{-}(x)\right)^{-1} G_{\varphi} \left(F^{+}(x)\right)^{-1} = \left(F^{-}(x)\right)^{-1} R_{\varphi} F R_{\varphi}^{-1} \left(F^{+}(x)\right)^{-1}.
$$

\noindent (6) There exist a small parameter $\varepsilon=\varepsilon(\varphi)$ (more exactly its value will be described later), such that for all $x\in {\mathbb R}$ and any finite $\varphi$
\begin{equation}
\label{constr2}
G_{1, \varphi}(x) = I + \varepsilon N_{\varphi}(x),
\end{equation}
and matrix $N_{\varphi}(x)$
is bounded and locally H\"older continuous on ${\mathbb R}$.


Note that by assumption each entry of the matrix $N_{\varphi}(x)$ has a limit when $|x| \rightarrow +\infty$, i.e. there exists the value $N_{\varphi}(\infty)$.
Note also that no commutativity of the involved matrices is assumed.

Let us look for factorization of the matrix $G_{1, \varphi}(x)$ in the form
\begin{equation}
\label{repr1}
G_{1, \varphi}(x) = I + \varepsilon N_{\varphi}(x) = \left(I + \varepsilon N_{1,\varphi}^{-}(x)\right) \left(I + \varepsilon N_{1,\varphi}^{+}(x)\right).
\end{equation}
Comparing terms at different powers of $\varepsilon$ we get, in particular,  the following relation for determination of factors
$N_{1,\varphi}^{-}$, $N_{1,\varphi}^{+}$:
\begin{equation}
\label{jump1}
N_{1,\varphi}^{-}(x) + N_{1,\varphi}^{+}(x) =  N_{\varphi}(x), \; x\in {\mathbb R}.
\end{equation}
It is customary to denote
$$
{M}_{0,\varphi}(x) \equiv N_{\varphi}(x).
$$
The jump boundary value problem (\ref{jump1}) has a solution represented in terms of a slight modification of the matrix-valued Cauchy type integral
(see, e.g., \cite{LitSpi87}, \cite{Mus68})
\begin{equation}
\label{sol_1}
N_{1,\varphi}^{\mp}(z) = \frac{1}{2} {M}_{0,\varphi}(\infty) \mp \frac{1}{2 \pi i} \int\limits_{-\infty}^{+\infty} \frac{{M}_{0,\varphi}(t) - {M}_{0,\varphi}(\infty)}{t - z} dt \equiv   \frac{1}{2} {M}_{0,\varphi}(\infty) \mp \left({\bf C}_{0}  {M}_{0,\varphi}\right)(z).
\end{equation}
The above modification is proposed in order to avoid extra discussion of the convergence of the above integrals. In this form the integrals are convergent automatically. Moreover, its boundary values $$\left({\bf C}_{0}  {M}_{0,\varphi}\right)^{\mp}(x) = \lim\limits_{{\mathrm{Im}} z \rightarrow \mp 0}  \left({\bf C}_{0}  {M}_{0,\varphi}\right)(z)$$
satisfy Sokhotsky-Plemelj formulas, i.e.
$$
\left({\bf C}_{0}  {M}_{0,\varphi}\right)^{\mp}(x) = \frac{1}{2} {M}_{0,\varphi}(x) - \frac{1}{2} {M}_{0,\varphi}(\infty) \mp \frac{1}{2 \pi i} \int\limits_{-\infty}^{+\infty} \frac{{M}_{0,\varphi}(t) - {M}_{0,\varphi}(\infty)}{t - x} dt =
$$
\begin{equation}
\label{sol_1_SP}
=
\frac{1}{2} {M}_{0,\varphi}(x) - \frac{1}{2} {M}_{0,\varphi}(\infty) \mp \frac{1}{2} \left({\bf S}_{0}  {M}_{0,\varphi}\right)(x),
\end{equation}
or
\begin{equation}
\label{sol_1_SP1}
N_{1,\varphi}^{\mp}(x) = \frac{1}{2} {M}_{0,\varphi}(x) \mp \frac{1}{2} \left({\bf S}_{0} {M}_{0,\varphi}\right)(x),
\end{equation}
where ${\bf S}_{0} $ is the singular integral operator along the real line with density ${M}_{0,\varphi}(t) - {M}_{0,\varphi}(\infty)$. It follows from \cite[n. 4.6]{Gak77} that both matrices $N_{1,\varphi}^{-}(x)$, $N_{1,\varphi}^{+}(x)$ satisfy H\"older conditions on $\overline{\mathbb R}$, are bounded there with $N_{1,\varphi}^{\mp}(\infty) = \frac{1}{2} {M}_{0,\varphi}(\infty)$, and possess an analytic continuation into lower $\Pi^{-}$ and upper $\Pi^{+}$ half-planes, respectively. Surely, its product is also  H\"older conditions on $\overline{\mathbb R}$ and bounded.

Let us refine the factorization of the matrix $G_{1, \varphi}(x)$, i.e. look for factorization in the form
\begin{equation}
\label{repr2}
G_{1, \varphi}(x) = I + \varepsilon N_{\varphi}(x) = \left(I + \varepsilon N_{1,\varphi}^{-}(x) + \varepsilon^2 N_{2,\varphi}^{-}(x)\right) \left(I + \varepsilon N_{1,\varphi}^{+}(x) + \varepsilon^2 N_{2,\varphi}^{-}(x)\right),
\end{equation}
where $N_{1,\varphi}^{-}(x)$, $N_{1,\varphi}^{+}(x)$ are those found at the previous step.

Comparing terms at different powers of $\varepsilon$ we get, in particular,  the following relation for determination of factors
$N_{2,\varphi}^{-}$, $N_{2,\varphi}^{+}$:
\begin{equation}
\label{repr2_1}
N_{1,\varphi}^{-}(x) N_{1,\varphi}^{+}(x) + N_{2,\varphi}^{-}(x) + N_{2,\varphi}^{+}(x) = 0.
\end{equation}
Denoting
$$
{M}_{1,\varphi}(x) \equiv - N_{1,\varphi}^{-}(x) N_{1,\varphi}^{+}(x),
$$
we arrive at the following jump boundary value problem
\begin{equation}
\label{jump2}
 N_{2,\varphi}^{-}(x) + N_{2,\varphi}^{+}(x) = {M}_{1,\varphi}(x),\; x\in {\mathbb R},
\end{equation}
with already known  right hand-side. Solution of this problem is given by the formula similar to (\ref{sol_1})
\begin{equation}
\label{sol_2}
N_{2,\varphi}^{\mp}(z) = \frac{1}{2} {M}_{1,\varphi}(\infty)
\mp \left({\bf C}_{0} {M}_{1,\varphi}\right)(z).
\end{equation}
It has the same properties as the solution of (\ref{sol_1}), in particular, its boundary values satisfy the relation
\begin{equation}
\label{sol_2_SP1}
N_{2,\varphi}^{\mp}(x) = \frac{1}{2} {M}_{1,\varphi}(x) \mp \frac{1}{2} \left({\bf S}_{0}  {M}_{1,\varphi}\right)(x),
\end{equation}
and $N_{2,\varphi}^{\mp}(\infty) = - \frac{1}{8} N_{\varphi}^{2}(\infty)$.

One can proceed in the same manner. Thus on the $k$-th step we use the representation
$$
G_{1, \varphi}(x) = I + \varepsilon N_{\varphi}(x) =
$$
\begin{equation}
\label{repr_k}
= \left(I + \varepsilon N_{1,\varphi}^{-}(x) + \ldots + \varepsilon^k N_{k,\varphi}^{-}(x)\right) \left(I + \varepsilon N_{1,\varphi}^{+}(x) + \ldots + \varepsilon^k N_{k,\varphi}^{-}(x)\right),
\end{equation}
where $N_{1,\varphi}^{-}(x), \ldots, N_{k-1,\varphi}^{-}(x)$, $N_{1,\varphi}^{+}(x), \ldots, N_{k-1,\varphi}^{-}(x)$ are found at the previous steps.
It leads to the jump boundary value problem
\begin{equation}
\label{jump_k}
 N_{k,\varphi}^{-}(x) + N_{k,\varphi}^{+}(x) = {M}_{k-1,\varphi}, \; x\in {\mathbb R},
\end{equation}
where
$$
{M}_{k-1,\varphi} = -\left[N_{1,\varphi}^{-}(x) N_{k-1,\varphi}^{+}(x) +
N_{2,\varphi}^{-}(x) N_{k-2,\varphi}^{+}(x) + \ldots + N_{k-1,\varphi}^{-}(x) N_{1,\varphi}^{+}(x)\right].
$$
 Solution of this problem is given by the formula similar to (\ref{sol_1})
(or  to (\ref{sol_2})).

Thus the factorization of the matrix function $G_{1, \varphi}(x)$ is given in the form of asymptotic series
\begin{equation}
\label{asy_fact}
G_{1, \varphi}(x) = \left({I}+\sum\limits_{k=1}^{\infty} \varepsilon^k N_{k,\varphi}^{-}(x) \right)
\left({I}+\sum\limits_{k=1}^{\infty} \varepsilon^k N_{k,\varphi}^{+}(x) \right),
\end{equation}
where the pair $N_{k,\varphi}^{-}(x)$, $N_{k,\varphi}^{+}(x)$ is the unique solution to the jump problem
(\ref{jump_k}) for any $k\in {\mathbb N}$.

The following Theorem gives conditions when this asymptotic factorization becomes  an explicit one, i.e. gives convergence conditions for the asymptotic series involved.

\begin{theorem}
\label{converge}
{Let the parameter $\varepsilon$ satisfies the inequality
\begin{equation}
\label{conv_cond}
|\varepsilon| \leq 1/A
\end{equation}
with the constant {$A=A(\varphi)$} being equal to
\begin{equation}
\label{conv_cond1}
A = \|N_{\varphi}(\cdot)\|_{\mu} (1 + C_{\mu})^2,
\end{equation}
$\|N_{\varphi}(\cdot)\|_{\mu}$ being the norm of the matrix function $N_{\varphi}(x)$ in the H\"older space $H_{\mu}$
equal to the maximum of the norms of its entries,
and $C_{\mu}$ being the norm of the singular integral operator ${\bf S}_{0} : H_{\mu} \rightarrow H_{\mu}$.}

{Then both series in the right-hand side of (\ref{asy_fact}) converge for all $x\in {\mathbb R}$.}
\end{theorem}
{$\triangleleft$
Boundedness of the modified singular integral operator ${\bf S}_{0}$ in H\"older spaces follows (see, e.g., \cite[p. 48]{Gak77}) from the boundedness of the ``standard'' singular integral operator
${\bf S}$ (Hilbert transform on the real line) in these spaces (the later is well-known, see  \cite{Ale75} for the exact value  of the  norm of $\|{\bf S}\|_{H_{\mu} \rightarrow H_{\mu}}$). Let us denote the norm of ${\bf S}_{0}$ in H\"older space $H_{\mu}({\mathbb R})$ by $C_{\mu}$, i.e.
$$
C_{\mu} = \|{\bf S}_{0}\|_{H_{\mu} \rightarrow H_{\mu}}.
$$
Then we have the following series of estimates
$$
\|N_{1,\varphi}^{\mp}(\cdot)\|_{\mu} \leq \alpha_1 \|N_{\varphi}(\cdot)\|_{\mu} (1 + C_{\mu}), \; {\mathrm{where}}\; \alpha_1 = \frac{1}{2},
$$
$$
\|N_{2,\varphi}^{\mp}(\cdot)\|_{\mu} \leq \frac{1}{2} \|M_{1,\varphi}(\cdot)\|_{\mu} (1 + C_{\mu}),
$$
and
$$
\|M_{1,\varphi}(\cdot)\|_{\mu} \leq \left(\alpha_1 \|N_{\varphi}(\cdot)\|_{\mu} (1 + C_{\mu})\right)^2,
$$
i.e.
$$
\|N_{2,\varphi}^{\mp}(\cdot)\|_{\mu} \leq \alpha_2 \|N_{\varphi}(\cdot)\|_{\mu}^{2} (1 + C_{\mu})^{3}, \; {\mathrm{where}}\; \alpha_2 = \frac{1}{2} \alpha_1^2.
$$
Finally, for each $k \geq 2$
$$
\|N_{k,\varphi}^{\mp}(\cdot)\|_{\mu} \leq \frac{1}{2} \|M_{k-1,\varphi}(\cdot)\|_{\mu} (1 + C_{\mu}),
$$
and
$$
\|M_{k-1,\varphi}(\cdot)\|_{\mu} \leq \left(\alpha_1 \alpha_{k-1} + \alpha_2 \alpha_{k-2} + \ldots + \alpha_{k-1} \alpha_1\right) \|N_{\varphi}(\cdot)\|_{\mu}^{k - 1} \left(1 + C_{\mu}\right)^{2k - 2},
$$
i.e.
$$
\|N_{k,\varphi}^{\mp}(\cdot)\|_{\mu} \leq \alpha_k \|N_{\varphi}(\cdot)\|_{\mu}^{k} (1 + C_{\mu})^{2k - 1}, \; {\mathrm{where}}\; \alpha_k = \frac{1}{2} (\alpha_1 \alpha_{k-1} + \ldots + \alpha_{k-1} \alpha_1).
$$
Few first coefficients $\alpha_k$ we can calculate explicitly, namely, $\alpha_1 = 1/2$, $\alpha_2 = 1/8$, $\alpha_3 = 1/16$. As for coefficients with large enough
indices we can proof by induction that
$$
\alpha_k < \frac{1}{16(k-3)}, \;\;\; \forall k \geq 12.
$$
Therefore
$$
\|\varepsilon^k N_{k,\varphi}^{\mp}(\cdot)\|_{\mu} \leq |\varepsilon|^k \frac{1}{32(k-3)(1 + C_{\mu})} \left(\|N_{\varphi}(\cdot)\|_{\mu} (1 + C_{\mu})\right)^{k}, \; \forall k \geq 12.
$$
Since the sequence $\sqrt[k]{\frac{1}{32(k-3)(1 + C_{\mu})}} \leq 1$ is increasing for sufficiently large $k$ and
$$
\lim\limits_{k \rightarrow\infty} \sqrt[k]{\frac{1}{32(k-3)(1 + C_{\mu})}} = 1,
$$
then the convergence of the series
$$
\left({I}+\sum\limits_{k=1}^{\infty} \varepsilon^k N_{k,\varphi}^{-}(x) \right),
\left({I}+\sum\limits_{k=1}^{\infty} \varepsilon^k N_{k,\varphi}^{+}(x) \right)
$$
for all $x\in {\mathbb R}$ follows from (\ref{conv_cond}).
$\triangleright$}

\begin{remark}
\label{rem_conv}
It follows from the standard properties of the Cauchy type integral and singular integral with Cauchy kernel
that conditions of Theorem \ref{converge} guarantee convergence of the series in the right-hand side of (\ref{asy_fact})
in the half-planes ${\Pi}^{-}$, ${\Pi}^{+}$, respectively.
\end{remark}

\begin{remark}
\label{infinity}
In fact, the decay of the second term in the right-hand side of (\ref{constr2}) at infinity
follows from the asymptotic relations (\ref{infinity_unit}) which in turn follows from the properties of matrices
of the considered class and the proposed construction).
\end{remark}

\begin{remark}
\label{smallness}
If the number $A = A(\varphi)$ in Theorem \ref{converge} is small enough, i.e.
\begin{equation}
\label{conv_cond1_without}
A = \|N_{\varphi}(\cdot)\|_{\mu} (1 + C_{\mu})^2<1,
\end{equation}
then the results remains valid for $\varepsilon = 1$ and the described procedure is working then without any changes.
\end{remark}


\subsection{Special case}
\label{special}

Let us consider the problem of factorization of
$2\times 2$ invertible matrices from a  subclass of ${\mathcal G} {K}_{2}$, namely
\begin{equation}
\label{matrix}
G_{\varphi}(x) = \left(
\begin{array}{cc}
p(x) & q(x) e^{i x \varphi} \\
q(x) e^{-i x \varphi} & p(x)
\end{array}
\right),
\end{equation}
 given on the real line ($x\in {\mathbb R}$) and depending on the real parameter $\varphi\in {\mathbb R}$.

We assume that the following assumptions hold.

1) the entries $p(x), q(x)$ are real-valued H\"older continuous functions on $\overline{\mathbb R}$, i.e. $p, q\in H_{\mu}(\overline{\mathbb R})$;

2) the combinations of the functions are positive:
\begin{equation}
\label{positive}
p(x) \pm q(x) > 0, x\in {\mathbb R};
\end{equation}

3) the following limits exist
\begin{equation}
\label{limit}
\lim\limits_{|x|\rightarrow +\infty} p(x) = 1,\;\;\; \lim\limits_{|x|\rightarrow +\infty} q(x) = 0;
\end{equation}

4) the following symmetry condition is valid
\begin{equation}
\label{symmetry}
\overline{G_{\varphi}(-x)} = G_{\varphi}(x);
\end{equation}

As in general case the factorization of matrices of type (\ref{matrix}) is motivated by certain problem of fracture mechanics. The considered matrices are similar
to those which are studied and explicitly factorized in \cite{Akt92}, \cite{FelGohKru04}, but certain conditions
of the above cited papers do not satisfy in our case.

Note, that even if one supposes that the  functions
$p(x), q(x)$ are meromorphically continued into semi-planes $\Pi^{-}$, $\Pi^{+}$, then it does not mean that these extended functions have finite number of zeroes and poles there.

\begin{remark}
\label{rem1}
Under conditions 1)--4) the matrix  (\ref{matrix}) admits the canonical factorization:
\begin{equation}
\label{fact_can}
G_{\varphi}(x) = G_{\varphi}^{-}(x) G_{\varphi}^{+}(x).
\end{equation}
\end{remark}
It follows, in particular, from condition 4) and \cite[p. 52]{Vor11}.

Let us start with factorization of an auxiliary matrix $F(x) = G_{0}(x)$
having no exponential term in their entries
\begin{equation}
\label{matrix_0}
F(x) = \left(
\begin{array}{cc}
p(x) & q(x) \\
q(x) & p(x)
\end{array}
\right).
\end{equation}


Note that
\begin{equation}
\label{matrix_01}
F(x) = P\left(
\begin{array}{cc}
p(x)+q(x) & 0 \\
0& p(x)-q(x)
\end{array}
\right)P,
\end{equation}
where the projector $P$ is defined
\[
P=\frac{1}{\sqrt{2}}\left(
\begin{array}{cc}
1 & 1 \\
1& -1
\end{array}
\right), 
\]


It follows from the conditions 1)--3) that both diagonal elements of the middle diagonal  matrix  have index equal to  zero. Hence they admit the representation
\begin{equation}
\label{bvp2_1}
\begin{array}{c}
p(x) + q(x) = \left(p(x) + q(x)\right)^{-} \cdot \left(p(x) + q(x)\right)^{+}; \\
p(x) - q(x) = \left(p(x) - q(x)\right)^{-} \cdot \left(p(x) - q(x)\right)^{+}
\end{array}
\end{equation}
with the factors of the form (see, e.g., \cite{Gak77})
\begin{equation}
\label{bvp2_2}
\begin{array}{ll}
 (p(z) + q(z))^{\mp} = \exp \left\{ \mp \frac{1}{2 \pi i} \int\limits_{-\infty}^{+\infty}
\frac{\log\, \left[p(\tau) + q(\tau)\right]}{\tau - z} d \tau\right\}, & z\in \Pi^{\mp}, \\
 (p(z) - q(z))^{\mp} = \exp \left\{ \mp \frac{1}{2 \pi i} \int\limits_{-\infty}^{+\infty}
\frac{\log\, \left[p(\tau) - q(\tau)\right]}{\tau - z} d \tau\right\}, & z\in \Pi^{\mp}
\end{array}
\end{equation}
and boundary values $ (p(x) + q(x))^{\mp}$ of the functions $ (p(z) + q(z))^{\mp}$ (and $ (p(x) - q(x))^{\mp}$ of the functions $ (p(z) - q(z))^{\mp}$ are determined by
using Sokhotsky-Plemelj formulas (see, e.g., \cite{Gak77}).

Therefore,  we obtain immediately the right  canonical factorization of the matrix $G_{0}(x)$:
\begin{equation}
\label{fact_can1}
F(x) = F^{-}(x) F^{+}(x),
\end{equation}
where
\begin{equation}
\label{fact_can_minus}
 F^{-}(x) = \frac{1}{\sqrt{2}}
\left( \begin{array}{lr}
 (p(x) + q(x))^{-} &  (p(x) - q(x))^{-} \\
 (p(x) + q(x))^{-} & -  (p(x) - q(x))^{-}
\end{array}
\right),
\end{equation}
\begin{equation}
\label{fact_can_plus}
 F^{+}(x) = \frac{1}{\sqrt{2}}
\left( \begin{array}{lr}
 (p(x) + q(x))^{+} &  (p(x) + q(x))^{+} \\
 (p(x) - q(x))^{+} & -  (p(x) - q(x))^{+}
\end{array}
\right).
\end{equation}
Now we can represent the initial matrix $G_{\varphi}(k)$ in the form:
\begin{equation}
\label{matrix1}
G_{\varphi}(x) = F^{-}(x) \left(F^{-}(x)\right)^{-1} G_{\varphi}(x) \left(F^{+}(x)\right)^{-1}
F^{+}(x) \equiv F^{-}(x) G_{1, \varphi}(x)  F^{-}(x)
\end{equation}
and proceed with factorization of the matrix $ G_{1, \varphi}(x)$.

The inverse matrices $\left(F^{-}(x)\right)^{-1}$, $\left(F^{+}(x)\right)^{-1}$ are equal, respectively
\begin{equation}
\label{matrix_in_minus}
\left(F^{-}(x)\right)^{-1} = \frac{1}{\Delta_{-}}
\left( \begin{array}{lr}
 (p(x) - q(x))^{-} &  (p(x) - q(x))^{-} \\
 (p(x) + q(x))^{-} & -  (p(x) + q(x))^{-}
\end{array}
\right),
\end{equation}
$$
\Delta_{-} = {\frac{2}{\sqrt{2}} } (p(x) + q(x))^{-}  (p(x) - q(x))^{-},
$$
\begin{equation}
\label{matrix_in_plus}
\left(F^{+}(x)\right)^{-1} = \frac{1}{\Delta_{+}}
\left( \begin{array}{lr}
 (p(x) - q(x))^{+} &  (p(x) + q(x))^{+} \\
 (p(x) - q(x))^{+} & -  (p(x) + q(x))^{+}
\end{array}
\right),
\end{equation}
$$
\Delta_{+} = {\frac{2}{\sqrt{2}}}  (p(x) + q(x))^{+}  (p(x) - q(x))^{+}.
$$
Hence (for shortness we omit further argument $x$ of functions {$p$} and $q$)
\begin{equation}
\label{matrix2}
 G_{1, \varphi}  =  \!\!\!\frac{1}{\Delta}
\!\!\!\left( \begin{array}{cc}
\!\!\! (p - q)(2 p + q (e^{i \varphi k} + e^{-i \varphi k})) &\!\!\!
(p - q)^{-} (p + q)^{+} q (e^{-i \varphi k} - e^{i \varphi k}) \\
\!\!\! (p - q)^{+} (p + q)^{-} q (e^{i \varphi k} - e^{-i \varphi k}) &\!\!\!
 (p + q) (2 p - q (e^{i \varphi k} - e^{-i \varphi k}))
\end{array}
\!\!\!\right),
\end{equation}
$$
\Delta = {2} (p^2  - q^2).
$$
It is not hard to see that $G_{1, \varphi}$ is a sum of the unit matrix and
the matrix which is ``small'' for appropriate choice of $\varphi$. Hence, following
Remark \ref{smallness}, we rewrite the right-hand side of (\ref{matrix2}) as the following sum
\begin{equation}
\label{matrix3}
 G_{1, \varphi}  =  I + N_{\varphi},
\end{equation}
where
\begin{equation}
\label{matrix_rest}
N_{\varphi} = \left( \begin{array}{lr}
-\frac{2 q \sin^2 \frac{\varphi x}{2}}{(p + q)} &
- \frac{i q (p - q)^{-} (p + q)^{+} \sin \varphi x}{(p^2 - q^2)} \\
\frac{i q (p - q)^{+} (p + q)^{-} \sin \varphi x}{(p^2 - q^2)} &
 \frac{2 q \sin^2 \frac{\varphi x}{2}}{(p - q)}
\end{array}
\right).
\end{equation}
Note that the latter matrix can be written as the sum of two diagonal matrices, namely
\begin{equation}
\label{matrix_rest0}
N_{\varphi} =  N_{1} +   N_{2}
\end{equation}
and both matrices $N_{1}$, $N_{2}$ do not depend on $\varphi$
\begin{equation}
\label{matrix_rest1}
N_{1} = 2 \sin^2 \frac{\varphi x}{2} \left( \begin{array}{cc}
- \frac{q}{(p + q)}   &  0
 \\
0 &
\frac{q}{(p - q)}
\end{array}
\right),
\end{equation}
\begin{equation}
\label{matrix_rest2}
N_{2} = i  \sin \varphi x \left( \begin{array}{cc}
0 &
- \frac{q (p - q)^{-} (p + q)^{+}}{(p^2 - q^2)}  \\
 \frac{q (p - q)^{+} (p + q)^{-}}{(p^2 - q^2)} &
0
\end{array}
\right).
\end{equation}
Let us denote
\begin{equation}
\label{small_p}
\varepsilon_1 = {\varepsilon_1(\varphi)} \equiv \max\limits_{x\in {\mathbb R}}  \left|q(x) \sin\, \frac{\varphi x}{2}\right|.
\end{equation}
\begin{lemma}
\label{epsilon}
The parameter $\varepsilon_1$ can be taken smaller than any positive number $\delta$ by an appropriate choice of $\varphi$.
\end{lemma}
$\triangleleft$ Indeed, taking into account
the condition (\ref{limit})$_2$ one concludes that there exists $x_\delta>0$ such that for any $|x|\ge x_\delta$
$$
|q(x) \sin\, \frac{\varphi x}{2}| \le |q(x)|\le \delta.
$$
On the other hand, since $xq(x)$ belong to the space $H_{\mu}[-x_\delta,x_\delta] \subset  C[-x_\delta,x_\delta]$
there exists a constant $q_h>0$ such that
$$
|xq(x)|\le q_h, \quad \mbox{for all}\quad |x|\le x_\delta.
$$
Finally this means
$$
\left|q(x) \sin\, \frac{\varphi x}{2}\right| \leq |x q(x)|  \frac{\varphi}{2}\le \frac{\varphi}{2}q_h.
$$
Choosing $\varphi=2\delta/q_h$ we finish the proof.$\triangleright$

From the structure of the matrix $N_{\varphi}$ from (\ref{matrix_rest}), it follows that
that parameter $A = A(\varphi)$ discussed in Theorem 3.1 is smaller than 1 for an appropriate choice of $\varphi$.
It guarantees applicability of the general procedure in this special case.



\section{An example }
\label{Example}

We present here an example of $2\times 2$ matrix function $G_{\varphi} \in {\mathcal G} {K}_{2}$ for which the above discussed factorization does not involve Cauchy type integration for the components and auxiliary matrices.

Let
\begin{equation}
\label{ex_matrix}
G_{\varphi}(x) = \left(
\begin{array}{cc}
\frac{x^2 + 10}{x^2 + 1} & \frac{6}{x^2 + 1} e^{i x \varphi} \\
\frac{6}{x^2 + 1} e^{-i x \varphi} & \frac{x^2 + 10}{x^2 + 1}
\end{array}
\right),\; x\in {\mathbb R}.
\end{equation}
It is a special case of the matrix functions discussed in Subsec. \ref{special}. Here $p(x) =  \frac{x^2 + 10}{x^2 + 1}$, $q(x) = \frac{6}{x^2 + 1}$ and an auxiliary matrix $F(x) = G_{0}(x)$ has the form
\begin{equation}
\label{ex_matrix_1}
G_{0}(x) = \left(
\begin{array}{cc}
\frac{x^2 + 10}{x^2 + 1} & \frac{6}{x^2 + 1}  \\
\frac{6}{x^2 + 1}  & \frac{x^2 + 10}{x^2 + 1}
\end{array}
\right),\; x\in {\mathbb R}.
\end{equation}
The above matrix (\ref{ex_matrix}) satisfies all conditions 1) -- 4) of Subsec. \ref{special}, moreover,
the functions $p, q$ are not only H\"older continuous on ${\mathbb R}$ but infinitely differentiable.

It follows from condition 2) that
$$
{\mathrm{Ind}}\, \det\, G_{\varphi}(x) = {\mathrm{Ind}}\, \det\, G_{0}(x) = 0,
$$
besides, due to condition 4) the matrix $G_{\varphi}(x)$ admits the canonical factorization (if exists).

Canonical factorization of the matrix $ G_{0}(x)$
$$
G_{0}(x) = G_{0}^{-}(x) G_{0}^{+}(x)
$$
can be found in an explicit form:
\begin{equation}
\label{ex_G_0-}
G_{0}^{-}(x) = \frac{1}{\sqrt{2}}\left(
\begin{array}{cc}
\frac{x - 4 i}{x - i} & \frac{x - 2 i}{x - i}  \\
\frac{x - 4 i}{x - i}  & - \frac{x - 2 i}{x - i}
\end{array}
\right),
\end{equation}
\begin{equation}
\label{ex_G_0+}
G_{0}^{+}(x) = \frac{1}{\sqrt{2}}\left(
\begin{array}{cc}
\frac{x + 4 i}{x + i} & \frac{x + 4 i}{x + i}  \\
\frac{x + 2 i}{x + i}  & - \frac{x + 2 i}{x + i}
\end{array}
\right).
\end{equation}
Hence, one can calculate an auxiliary matrix
$$G_{1, \varphi}(x) = \left(G_{0}^{-}(x)\right)^{-1} G_{\varphi}(x) \left(G_{0}^{+}(x)\right)^{-1}
$$
in the following form
\begin{equation}
\label{ex_matrix_G_1}
G_{1, \varphi}(x) = \frac{1}{2 (x^2 + 16)(x^2 + 4)}\left(
\begin{array}{cc}
g_{1,1}(x) & g_{1,2}(x)  \\
g_{2,1}(x)  & g_{2,2}(x)
\end{array}
\right),
\end{equation}
where
\begin{equation}
\label{ex_g_11}
g_{1,1}(x) = 2(x^2 + 4)(x^2 + 10) + 6 (x^2 + 4) \left(e^{i \varphi x} + e^{-i \varphi x}\right),
\end{equation}
\begin{equation}
\label{ex_g_12}
g_{1,2}(x) = 12(x - 2 i)(x + 4 i) \left(-e^{i \varphi x} + e^{-i \varphi x}\right),
\end{equation}
\begin{equation}
\label{ex_g_21}
g_{2,1}(x) = 12(x + 2 i)(x - 4 i) \left(e^{i \varphi x} - e^{-i \varphi x}\right),
\end{equation}
\begin{equation}
\label{ex_g_22}
g_{2,2}(x) = 2(x^2 + 16)(x^2 + 10) + 6 (x^2 + 16) \left(-e^{i \varphi x} - e^{-i \varphi x}\right).
\end{equation}

It leads to the following representation of the matrix  
\begin{equation}
\label{ex_matrix_G_1_phi}
G_{1, \varphi}(x) = I +  N_1 + N_2,
\end{equation}
where we can take $\varepsilon = 1$ due to Theorem \ref{converge} (see also Remark \ref{smallness}), and
$$
N_1 = 6 i \sin^2 \frac{\varphi x}{2} \left(
\begin{array}{cc}
\frac{1}{x^2 + 16} & 0  \\
0  & - \frac{1}{x^2 + 4}
\end{array}
\right),
$$
$$
N_2 =  12 i \sin\, {\varphi x} \left(
\begin{array}{cc}
0 & - \frac{1}{(x + 2 i)(x - 4 i)}\\
\frac{1}{(x - 2 i)(x + 4 i)}  &  0
\end{array}
\right).
$$
This follows from the consideration presented in the previous Section. 

Following general scheme of Subsec. \ref{general} we factorize $G_{1, \varphi}(x)$ at first in the form (\ref{jump1})
where the first terms of factorization should be computed by the formula (\ref{sol_1}).
Following the asymptotic procedure described above, it is sufficient
at the first stage to factorize the matrix
\begin{equation}
\label{ex_fact1_2}
M_{0, \varphi}(x) = 6 i \left(
\begin{array}{cc}
\frac{\sin^2 \frac{\varphi x}{2}}{x^2 + 16} & - \frac{2 \sin\, {\varphi x}}{(x + 2 i)(x - 4 i)}  \\
 \frac{2 \sin\, {\varphi x}}{(x - 2 i)(x + 4 i)} & - \frac{\sin^2 \frac{\varphi x}{2}}{x^2 + 4}
\end{array}
\right).
\end{equation}

For this particular case, instead of using the Cauchy integrals, one can factorize each entry
$n_{ij}(x), i = 1, 2,$ of the matrix by using decomposition in simple fraction
and Taylor formula.
Combining the obtained results we get the following representation of matrix $ N_{1,\varphi}^{-}(x)$, $ N_{1,\varphi}^{+}(x)$ (first components of asymptotic factorization)
\begin{equation}
\label{ex_fact1-un}
N_{1,\varphi}^{-}(x) =  6 i \left(
\begin{array}{cc}
\frac{- e^{- i \varphi x} - e^{- 4 \varphi} + 2}{32 i (x - 4 i)} + \frac{e^{- i \varphi x} - e^{- 4 \varphi }}{32 i (x + 4 i)} & \frac{- i e^{- i \varphi x} + i e^{- 4 \varphi}}{6 i (x - 4 i)}
+ \frac{ i e^{- i \varphi x} - i e^{- 2 \varphi}}{6 i (x + 2 i)}\\[2mm]
\frac{i e^{- i \varphi x} - i e^{- 2 \varphi}}{6 i (x - 2 i)}
- \frac{ i e^{- i \varphi x} - i e^{- 4 \varphi}}{6 i (x + 4 i)} & \frac{e^{- i \varphi x} + e^{- 2 \varphi} - 2}{16 i (x - 2 i)} - \frac{e^{- i \varphi x} - e^{- 2 \varphi }}{16 i (x + 2 i)}
\end{array}
\right).
\end{equation}
\begin{equation}
\label{ex_fact1+un}
N_{1,\varphi}^{+}(x) =  6 i \left(
\begin{array}{cc}
\frac{e^{ i \varphi x} + e^{- 4 \varphi} - 2}{32 i (x + 4 i)} - \frac{e^{i \varphi x} - e^{- 4 \varphi }}{32 i (x - 4 i)} & \frac{- i e^{ i \varphi x} + i e^{- 2 \varphi}}{6 i (x + 2 i)}
+ \frac{ i e^{ i \varphi x} - i e^{- 4 \varphi}}{6 i (x - 4 i)}\\[2mm]
\frac{i e^{i \varphi x} - i e^{- 4 \varphi}}{6 i (x + 4 i)}
- \frac{ i e^{ i \varphi x} - i e^{- 2 \varphi}}{6 i (x - 2 i)} & \frac{- e^{ i \varphi x} - e^{- 2 \varphi} + 2}{16 i (x + 2 i)} + \frac{e^{i \varphi x} - e^{- 2 \varphi }}{16 i (x - 2 i)}
\end{array}
\right).
\end{equation}
Note that both plus ($N_{1,\varphi}^{+}$) and minus ($N_{1,\varphi}^{+}$) matrix-functions vanish at infinity in the corresponding half-plane.
 The procedure can be performed further as it was described in the general case.

\section{Outlook and discussions }

To illuminate the efficiency of the proposed procedure, we present here the numerical results related to the previous example (Sec. \ref{Example})
showing the quality of the factorization if one decides to restrict the approximation to the first asymptotic term only.

\begin{figure}[h!]
  \begin{center}
    \includegraphics [scale=0.35]{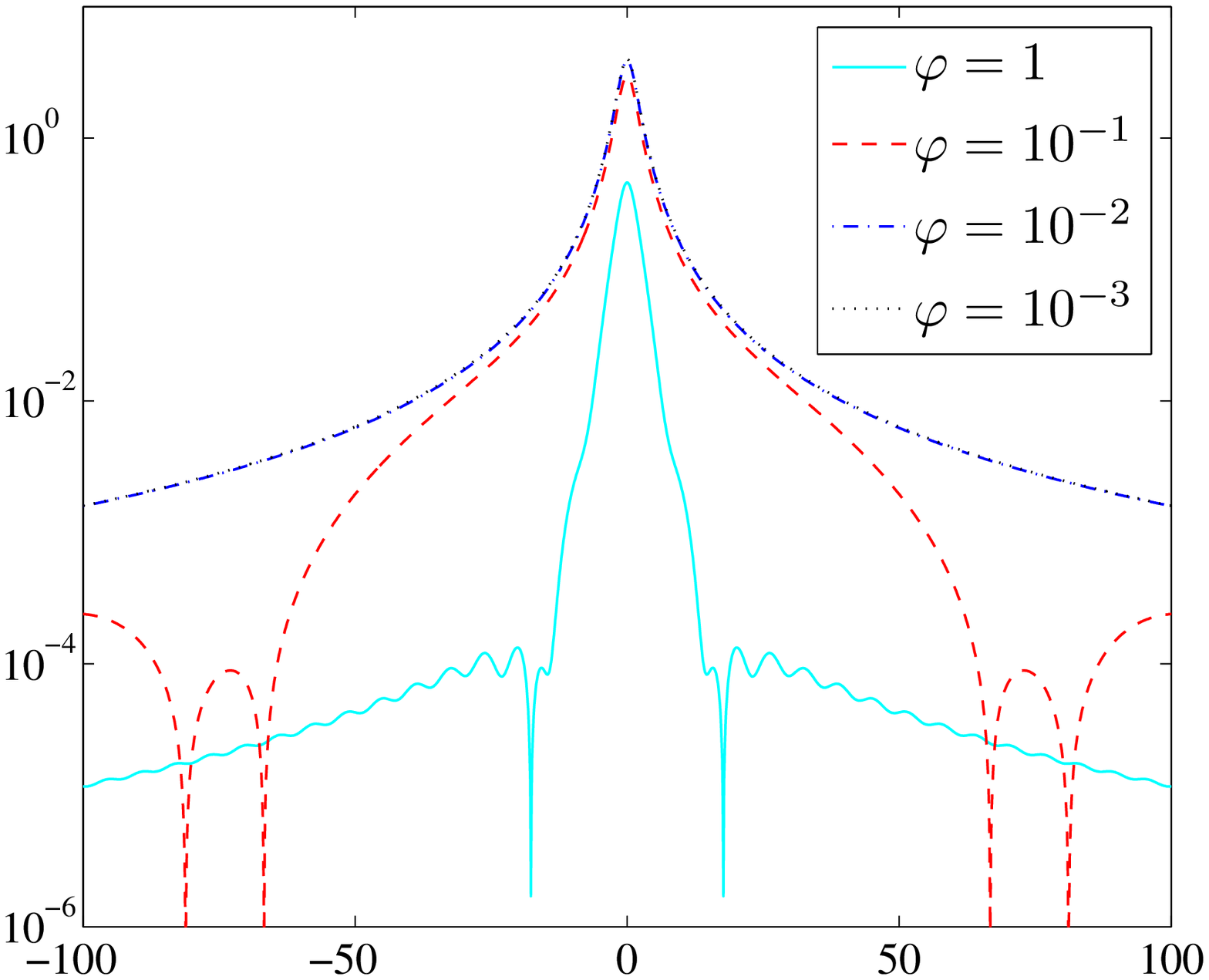}\includegraphics [scale=0.35]{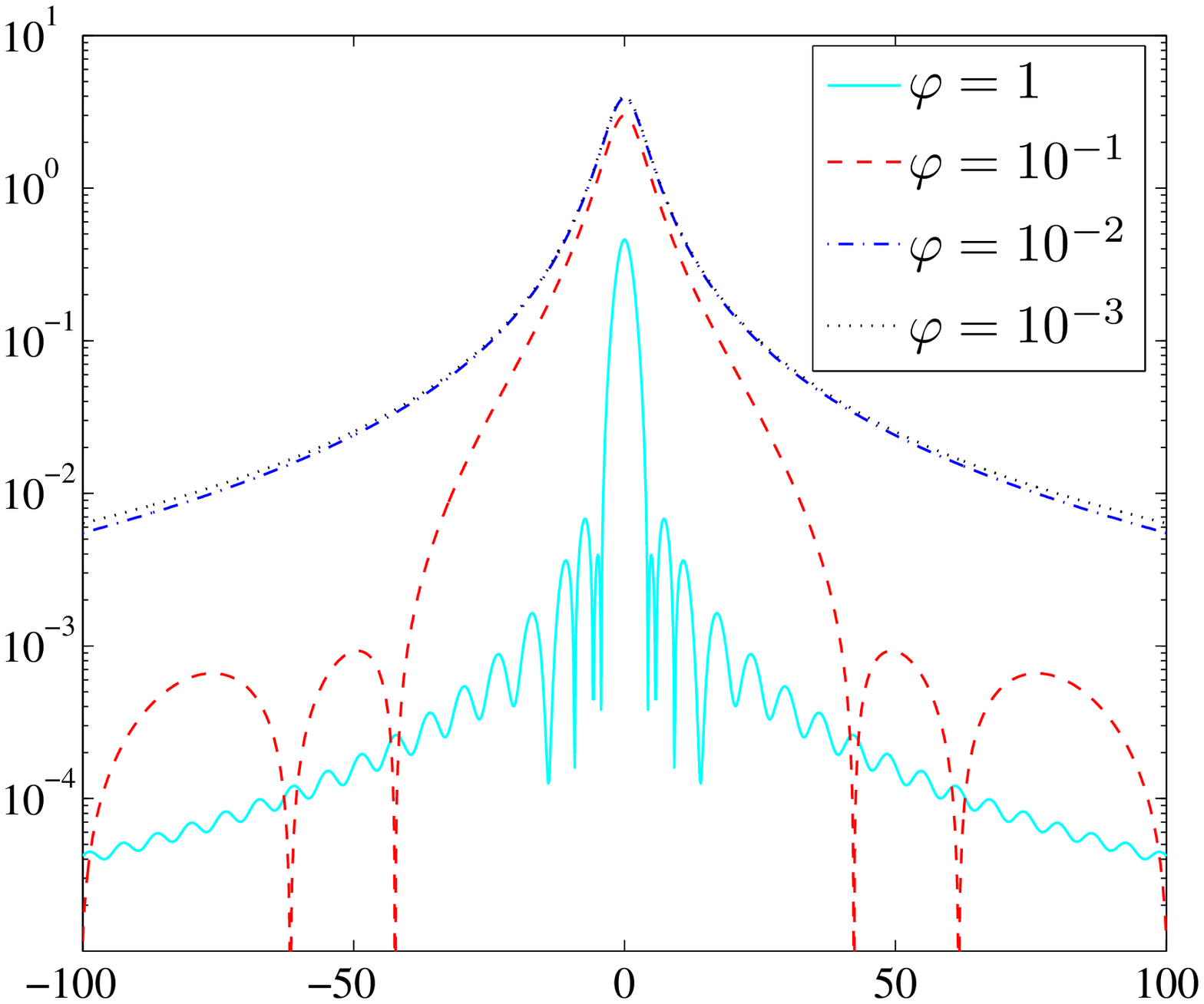}
       \put(-100,-5){\small$x$}
       \put(-365,120){$\varphi^{-2}\Delta k_{11}$}
       \put(-165,120){$\varphi^{-2}\Delta k_{22}$}
       \put(-300,-5){\small$x$}

  \end{center}
    \caption{Absolute error in the diagonal elements appeared by replacing the exact factorization with the only first asymptotic terms
    (\ref{ex_fact1-un})  and (\ref{ex_fact1+un}) for various values of the parameter $\varphi$.}
\label{f1}
\end{figure}

\begin{figure}[h!]
  \begin{center}
    \includegraphics [scale=0.35]{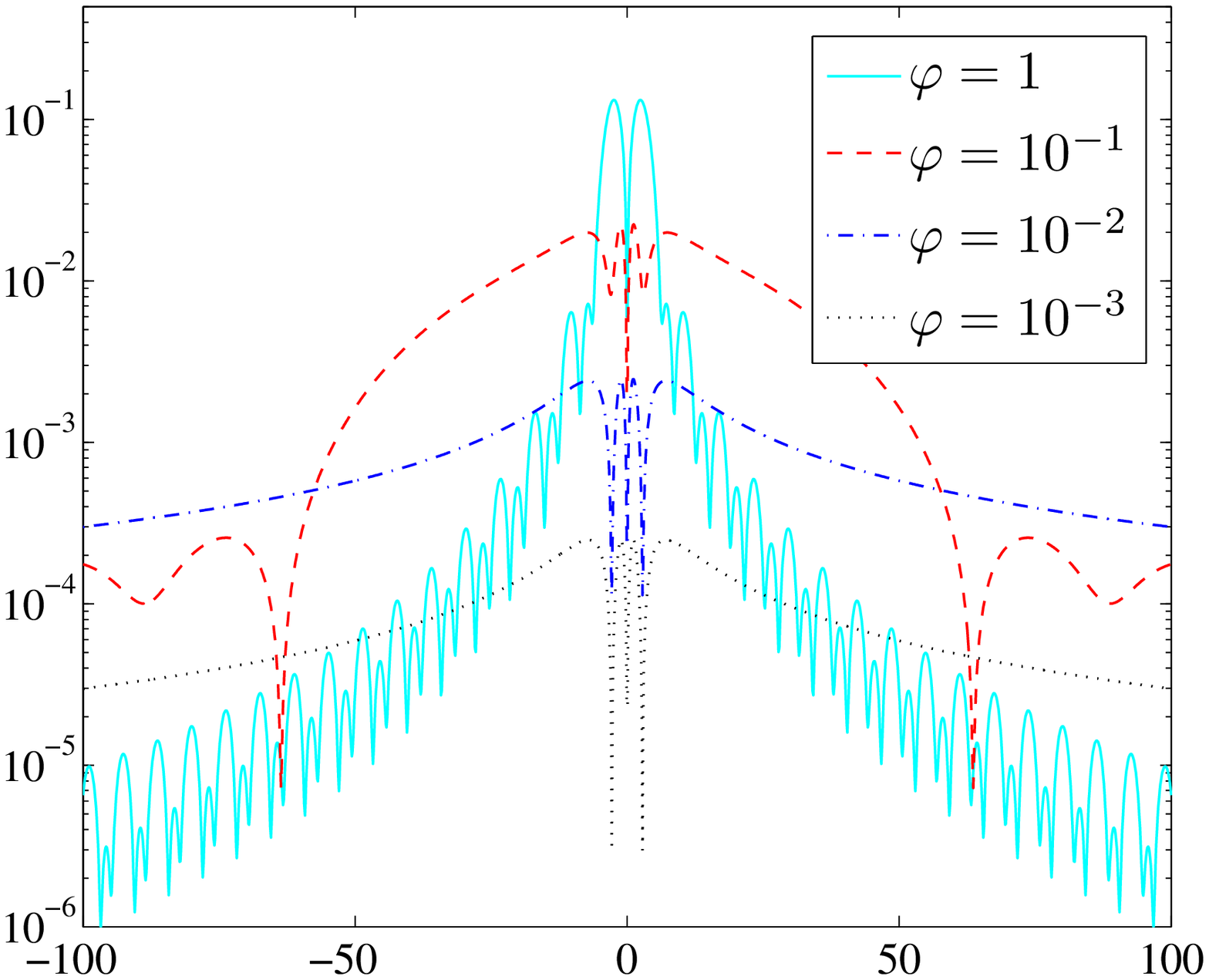}
    \begin{picture}(0,0)(0,0)
       \put(-100,-5){\small$x$}
         \put(-165,120){$\varphi^{-2}\Delta k_{ij}$}
\end{picture}

  \end{center}
    \caption{Absolute error in the off-diagonal elements appeared by replacing the exact factorization with the only first terms (\ref{ex_fact1-un})  and (\ref{ex_fact1+un}) for various values of the parameter $\varphi$.}
\label{f2}
\end{figure}

In Fig.~\ref{f1} and Fig.~\ref{f2} we present the normalized absolute error of the factorization using the only first asymptotic term
for different values of the parameter $\varphi$. We compute the errors for each component of reminder, $\Delta K$, that is the difference between the exact factorization and its asymptotic approximation along the real axis
\[
\Delta K(x)= G_{1, \varphi}(x)-\left( I + N_{1,\varphi}^{-}\right)  \left( I + N_{1,\varphi}^{+}\right)=N_{1,\varphi}^{-} N_{1,\varphi}^{+}.
\]
As it follows from the properties of the Cauchy type integral \cite{Mus68}, the obtained estimates are valid also into the upper and lower half-planes.

In Fig.~\ref{f1} we depict the error related to the diagonal elements of the reminder while in Fig.~\ref{f2} we show the result for the off-diagonal element. Note, that since they are complex conjugate with respect to each other, it is enough to discuss only one of them. Unexpectedly, even for rather large value of the parameter $\varphi=1$, the error appeared to be not high,
while for smaller magnitudes of $\varphi$ it decays fast with argument and its larger value is concentrated only near the center of the coordinate.  Moreover, one can observe that for $j=1,2$
\[
\Delta k_{jj}(0)=4\varphi^2+O(\varphi^3),\quad \varphi\to 0.
\]

This result can be verified also analytically.

With decrease of the parameter $\varphi$, the amplitude of the error oscillations increases but their support moves out from the coordinate center to infinity. Thus, such oscillations
plays minor role and the numerical computation of the Cauchy integrals in the described procedure would not affect the accuracy of the computations.

It is not correct to say, however, that $\varphi$ is an optimal small parameter in the problem under consideration.
As it follows from Theorem \ref{converge}, not only the parameter $\varphi$ but the decay of the corresponding function (e.g., the function $q(x)$ in the special case) play an important role in the analysis and it is rather $\delta$ should be taken as the appropriate small
parameter.

To clarify this, let us note that two terms $N_1$ and $N_2$
are of different orders with respect to the small parameter $\varphi$. Indeed:
 $N_1(\varphi)=O(\varphi^2)$ and $N_2(\varphi)=O(\varphi)$, as $\varphi\to0$.
As a result, one could construct another first order approximation of the factorization basing only on the term $N_2$ instead of $N_1+N_2$.
This gives the same first order estimate in terms of the small parameter $\varphi$. Simple calculations give the following (first-order) factorization terms in this case:
\begin{equation}
\label{N_2_minus}
N_{1,\varphi}^{*-}(x) = \left(
\begin{array}{cc}
0 & -  \frac{i \left(e^{-i \varphi x} - e^{- 4 \varphi}\right)}{x - 4 i} + \frac{i \left(e^{-i \varphi x} - e^{- 2 \varphi}\right)}{x + 2 i} \\
\frac{i \left(e^{-i \varphi x} - e^{- 2 \varphi}\right)}{x - 2 i} - \frac{i \left(e^{-i \varphi x} - e^{- 4 \varphi}\right)}{x + 4 i} & 0
\end{array}
\right),
\end{equation}
\begin{equation}
\label{N_2_plus}
N_{1,\varphi}^{*+}(x) = \left(
\begin{array}{cc}
0 & -  \frac{i \left(e^{i \varphi x} - e^{- 2 \varphi}\right)}{x + 2 i} + \frac{i \left(e^{i \varphi x} - e^{- 4 \varphi}\right)}{x - 4 i} \\
\frac{i \left(e^{i \varphi x} - e^{- 4 \varphi}\right)}{x + 4 i} - \frac{i \left(e^{i \varphi x} - e^{- 2 \varphi}\right)}{x - 2 i} & 0
\end{array}
\right).
\end{equation}

We estimate then a quality of the approximation of the factorization of the matrix-function $G_{1, \varphi}(x)$
basing now solely on the first order term with respect to the parameter $\varphi$. The new reminder is
\[
\Delta K^*(x)= G_{1, \varphi}(x)-\left( I + N_{1,\varphi}^{*-}\right)  \left( I + N_{1,\varphi}^{*+}\right)=N_1(x)-N_{1,\varphi}^{*-}N_{1,\varphi}^{*+}.
\]
Note that the off-diagonal terms gives exact (identity) result for the terms (\ref{N_2_minus}), (\ref{N_2_plus}) in opposite to the factorization provided by
the terms (\ref{ex_fact1-un}), (\ref{ex_fact1+un}). Thus for the first glance, the latter approach looks less beneficia than the former.

\begin{figure}[h!]
  \begin{center}
    \includegraphics [scale=0.35]{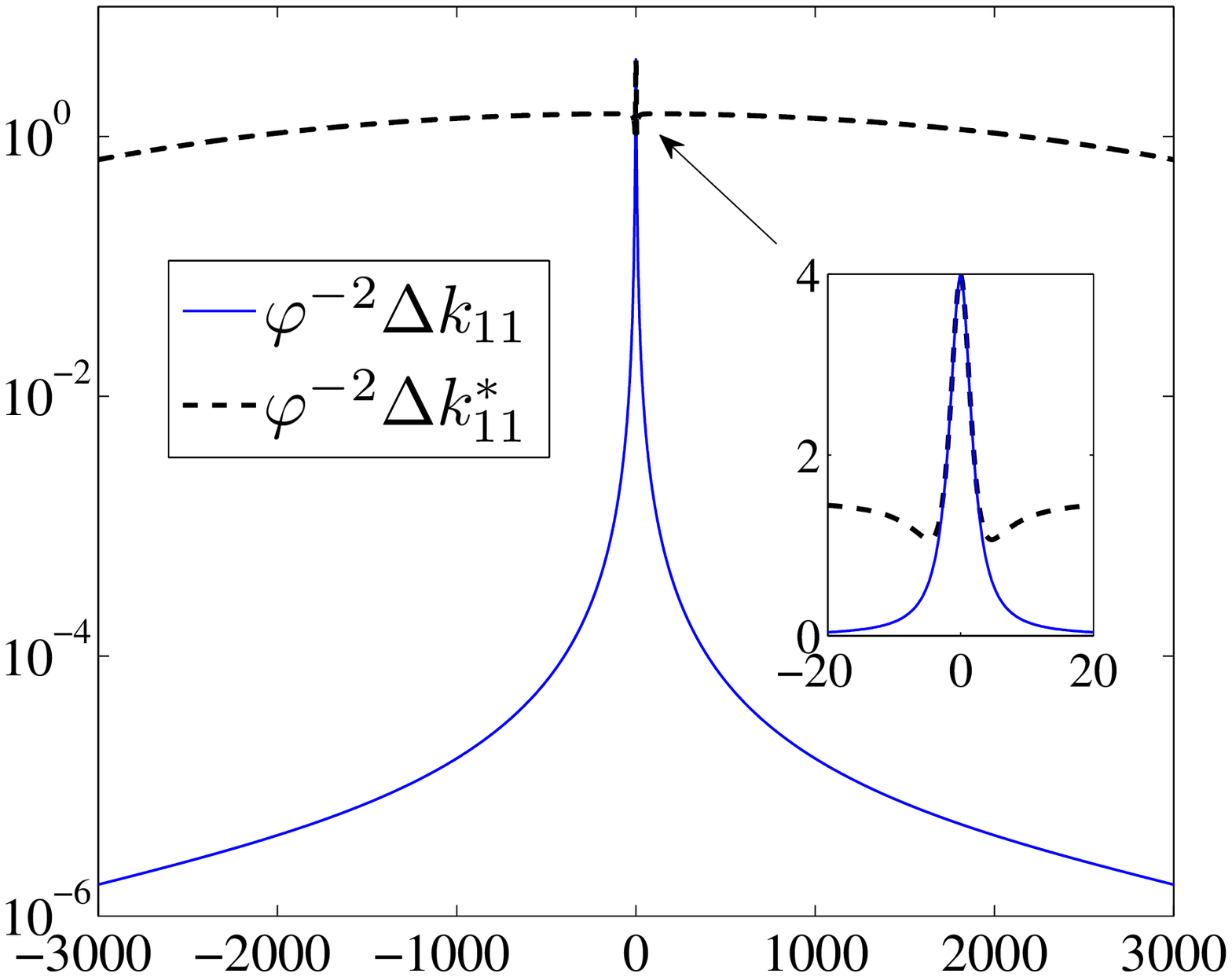}\includegraphics [scale=0.35]{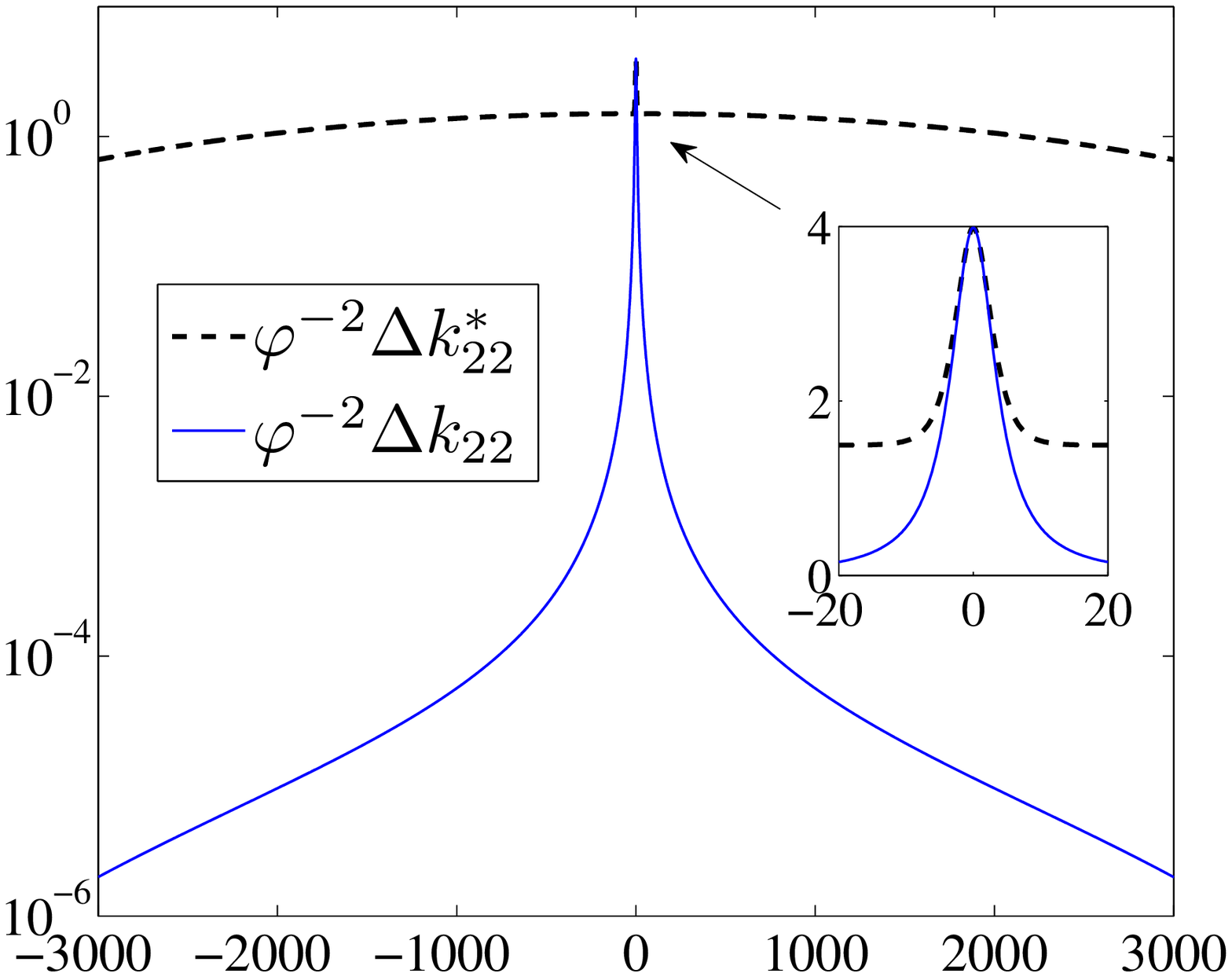}
       \put(-100,-5){\small$x$}
       \put(-300,-5){\small$x$}
  \end{center}
    \caption{Comparison of the errors in the diagonal elements and different approximations for the value of the parameter $\varphi=10^{-3}$.
    The errors corresponding to the factorization consisted of the
matrix-functions (\ref{N_2_minus}), (\ref{N_2_plus}) are given by doted line while those related to (\ref{ex_fact1-un}), (\ref{ex_fact1+un}) are depicted by the solid line.}
\label{f3}
\end{figure}

In Fig.~\ref{f3} we compare the errors related to the diagonal elements of the reminders for the value of the parameter $\varphi=10^{-3}$. The errors corresponding
to the
matrix-functions (\ref{N_2_minus}), (\ref{N_2_plus}) are given by doted line while those related to (\ref{ex_fact1-un}), (\ref{ex_fact1+un}) are depicted by the solid line.
One can observe a striking difference in the accuracy of those two approximations. While for a small values of the variable $|x|<1$ they are identical, for larger value of the argument the approximation given by the general procedure
provide much better accuracy than that based on the only first order term, $N_2$, with respect to the small parameter $\varphi$.
Moreover, if one decide to continue asymptotic expansion, it may become a real issue in numerical computations of the next asymptotic terms as the decay is very slow. Thus, the procedure suggested here is in a sense optimal.

Note, that proposed procedure is working in a general setting discussed in Sec. \ref{general} not only in the case of the Subsec. \ref{special}. In particular, the procedure is definitely valid
when all entries of the matrices
are quasi-polynimials, i.e. sum of different exponentials with meromorphic coefficients, providing that all conditions of the class ${\mathcal G} K_n$ are satisfied.

In the case of non-canonical factorization the algorithm becomes more cumbersome. This case will be considered elsewhere.

\vspace{3mm}
{\bf Acknowledgement.} {The work is supported by the FP7 PEOPLE IAPP project
N 251475 "HYDROFRAC".}


\end{document}